\documentclass[12pt]{article}

\usepackage[dvipsnames, usenames]{color}   
\usepackage{amsmath}
\usepackage{amssymb}
\begin{document}

\newcommand{\ff}{{\mathbb{ F\!}}}
\newcommand{\lff}{{\mathbb{ F\!}}^{\,\prime}}

\newtheorem{theorem}{Theorem}[section] 
\newtheorem{lemma}[theorem]{Lemma}     
\newtheorem{corollary}[theorem]{Corollary}
\newtheorem{proposition}[theorem]{Proposition}
\newtheorem{definition}[theorem]{Definition}
\newtheorem{defn}{Definition}
\newtheorem{example}[theorem]{Example} 
\newtheorem{remark}[theorem]{Remark}
\newtheorem{application}[theorem]{Application}
\newenvironment{proof}{\par\noindent {\bf Proof: }}{\hfill $\Box$ \\}

\renewcommand{\k}{K^{\times}}
\newcommand{\F}{\mathbb{F}}
\newcommand{\Fq}{\mathbb{F}_q}
\newcommand{\Fl}{\mathbb{F}_{\ell}}

\newcommand{\A}{\mathbb{A}}
\renewcommand{\P}{\mathbb{P}}
\newcommand{\leg}[1]{\left(\frac{#1}{\ell}\right)}

\newcommand{\di}{\mathrm{diag}}    
\newcommand{\ep}{\epsilon}         
\newcommand{\psl}{\mathrm{PSL}_2(q)}
\newcommand{\psll}{\mathrm{PSL}_2(\ell)}
\newcommand{\PSL}{\mathrm{PSL}_2(q^2)}

\newcommand{\SLG}{\mathrm{SL}_2(q^2)}
\newcommand{\slg}{\mathrm{SL}_2(\ell)}   
\newcommand{\slgl}{\mathrm{SL}_2(\ell)}
\newcommand{\ovr}{\overline}
\newcommand{\ovth}{\ovr{\theta}} 
\renewcommand{\eth}{E_{\theta}(x)}
\newcommand{\fth}{F_{\theta}(x)}

\newcommand{\vv} {{\boldsymbol{\nu}}}
\newcommand{\mm}{{\boldsymbol{\mu}}}
\newcommand{\la}{{\boldsymbol{\lambda}}}

\newcommand{\dd} {{\boldsymbol{\delta}}}
\renewcommand{\aa}{{\mathbf{a}}}
\newcommand{\xx}{{\bf x}}
\newcommand{\yy}{{\bf y}} 
\newcommand{\bl} {{\boldsymbol{\ell}}}
\newcommand{\ttt}{{\bf t}}
\newcommand{\sst}{{\mathcal S}}
\newcommand{\RR}{{\mathcal R}}
\newcommand{\N}{{\mathcal N}}
\newcommand{\T}{{\mathcal T}}
\newcommand{\nn}{{\vv}}
\newcommand{\sss}{{\bf s}}
\newcommand{\f}{{\mathbf{f}}}
\newcommand{\bbF}{\F}
\newcommand{\Fp}{{\mathbb{F\!}_p}}
\newcommand{\Z}{\mathbb{Z}}
\newcommand{\Q}{\mathbb{Q}}
\newcommand{\x}{\mathbf{x}}
\newcommand{\e}{\mathbf{e}}
\newcommand{\be}{\begin{equation}}
\newcommand{\eeq}{\end{equation}}
\newcommand{\bx}{{\mathbf{x}}}
\newcommand{\lcm}{\operatorname{lcm}}
\newcommand{\bb}{{\bf b}}
\newcommand{\ain}{\aa \in {\Fp}^n}
\newcommand{\bin}{\bb \in {\Fp}^n}
\newcommand{\CC}{\mathcal C}
\newcommand{\GG}{\mathcal G}

\newcommand{\Ssystem}{\mathcal{S}}
\newcommand{\Ematrix}{{E}}

\newcommand{\GL}{{\mathrm{GL}}_k(q)}
\newcommand{\GLF}{{\mathrm{GL}}_k(F)}
\newcommand{\qed}{{\hfill$\Box$}}
\newcommand{\pr}{\noindent{\bf Proof. \ }}
\newcommand{\lra}{\rightarrow}
\def\r{\rm res}
\def\l{{\cal l}}
\def\x{{\cal x}}
\def\y{{\cal y}}
\def\F{{\cal F}}
\def\C{{\cal C}}
\def\D{{\cal D}}
\def\L{{\cal L}}
\def\m{{\cal m}}
\def\P{{\cal P}}
\def\N{{\cal N}}
\def\O{{\cal O}}
\def\G{{\cal G}}
\def\Z{{\cal Z}}
\def\X{{\cal X}}
\def\M{{\cal M}}
\def\Y{{\cal Y}}
\def\x{{\bf x}}
\def\y{{\bf y}}
\def\u{{\bf u}}
\def\v{{\bf v}}
\def\w{{\bf w}}
\def\h{{\bf h}}
\def\c{{\bf c}}

\def\s{{\rm Sch}}
\def\cl{{\bf{\rm Cl}}(X)}
\def\div{{\bf{\rm Div}}(X)}

\def\closK{{\overline K}}
\def\Gal{{\rm Gal}}
\def\Div{{\rm Div}}
\def\DivK{{\rm Div}_K}
\def\deg{{\rm deg}}
\def\dg{C_{\Omega}}
\def\rsdg{C_{\L}}

\unitlength = 0.5cm
\def\map#1#2{\buildrel#1\over#2}
\def\N{I\!\!N}
\def\F{I\!\!F}
\def\Z{Z\!\!\!Z}
\def\Q{I\!\!\!Q}
\def\R{I\!\!R}
\def\C{I\!\!\!\!C}
\def\H{I\!\!H}

\author{H. F. Mattson, jr.}
\title{ New 5-Designs---revisited } 
\maketitle

The main purpose of this note is to clarify some portions of,
and correct some errors in,  
{\it New 5-Designs}, \cite{NFD}, hereafter ``{\sc nfd}'':  
The definitions of the
quadratic-residue codes in {\sc nfd} and their relationships to 
each other are clarified, as is their extension. The 
proof of the Gleason-Prange theorem is here made clearer.  

Determining the sign in the extension to the ``infinite''
coordinate but left ambiguous in {\sc nfd} led to this note.
  
Some examples calculated in a current work \cite{JM} applying {\sc nfd} 
to modular representations of groups made it desirable to determine the 
above-mentioned sign.

Two other matters first.

\medskip

{\sc Chebotar\"{e}v's theorem.}
One of the results in {\sc nfd}, pp. 127-128, was a new proof of the widely 
known Chebotar\"{e}v's 
theorem on the roots of unity.  At the time, the 1960s, we had no idea 
it was known, 
much less that it had its own name.  (It says: for a prime $\ell$ and 
a primitive $\ell^{th}$
root of unity $z$ over $\Q$, every subdeterminant of the $\ell \times \ell$
matrix $(z^{ij})$ is nonvanishing.)  But we
suggest our coding-theoretic proof is worth considering.  

We used this result to prove ``optimality''\footnote{A word long since displaced by the abbreviation ``MDS'', for a crackjaw term best left unsaid.} of all
cyclic $[\ell, k]$ codes
over $\F_{p^i}$ for almost all primes $p$.  (See Theorem 2.2 in {\sc nfd}.)

An article on
Chebotar\"{e}v's life and works, mentioning 
other proofs of his result, appears in \cite{cheb}.

\medskip

{\sc Orbits.}
An example at the end of {\sc nfd} has an error.  The matter is the action
of $\mathrm{PSL}_2(47)$ on the binary [48,24,12] code, in particular on the
codewords of minimum weight.  These codewords are in three orbits, but not
the orbits stated in {\sc nfd}.  A correct account appears in \cite{R1972}, 
pp. I-26{\it ff}.

\bigskip

{\sc The main thing}.  
 
The reader is assumed to be familiar with proof of the Gleason-Prange theorem
in {\sc nfd}.

Throughout we denote by $\ell$ an odd prime.  When we speak of a polynomial,
say, $a+bx+cx^2$,  as an element of a code, we understand that as an $n$-tuple 
it is $(a, b, c, 0, \ldots , 0)$.  The {\bf reverse} $f^*(x)$ of the non-0
polynomial  $f(x)$ of degree $d$  is defined as  
\[  f^*(x) := x^df(1/x).  \]
 
We begin with a result on cyclic codes.  We will freely represent codewords
either as $n$-tuples over a field or as polynomials.  

The setup:  $L$ is a field, and 
$z$ an  $n^{th}$ root of unity over $L$.  $K := L(z)$.

\bigskip

{\bf Lemma}. {\it Let $\varphi$ be any nontrivial linear functional from $K$ to $L$.  Let $h(x) \in L[x]$ of degree $k$ be the monic irreducible polynomial over $L$ with $z$ as root.  Let $C$  denote the cyclic code}
\[ C := \{F_z(c):=(\varphi(c), \varphi(cz), \ldots , \varphi(cz^{n-1}); c \in K\}.  \]
{\it Then $C$ is an $[n,k]$ cyclic code over $L$, and $C^{\perp} = (h(x))$, the $[n,n-k]$ cyclic code generated by $h(x)$.  And }
\[ C = \left(\frac{x^n -1}{h^*(x)}\right).  \]

\pr Let $h(x) = h_0 + h_1x + \cdots + h_kx^k$.  It is obvious that
$v := (h_0, h_1, \ldots ,h_k, 0, \ldots ,0)$ is orthogonal to $F_z(c)$.
The same is true of every cyclic shift of $v$.  As a polynomial $v$ is $h(x)$.
Thus $ (h(x)) \subseteq C^{\perp}  $.  That is, dim$(C) \leq k$.

For the reverse inequality, note that  $c \longrightarrow F_{z}(c)$   is a 
linear map from  $K$ to $L^n$.  If $c \neq 0$ then $\{c, cz, \ldots, cz^{n-1}\}$
spans $K/L$, so $F_z(c) \neq 0$.  Since $[K:L] = k$, that is the dimension 
of $C$.  Thus $C^{\perp} = (h(x))$.

The final assertion of the Lemma follows from the well known result that
when $x^n - 1 = f(x)g(x)$, the orthogonal code of the cyclic code $(g(x))$
is $(f^*(x))$.                \qed

\bigskip

Now to {\sc nfd}.  The purpose of this section is to clarify the definition
of the extended quadratic-residue codes and the proof of the Gleason-
Prange theorem.  Unexplained notations are taken from {\sc nfd}.

Let us clean the slate by permanently using the factoring  
\[  x^{\ell} - 1 = (x-1)f(x)g(x),    \]
instead of the $g_1(x)g_2(x)$ in {\sc nfd}, page 129.  Now $z$ is a primitive
$\ell^{th}$ root of 1 over $\Q$, and all takes place in  $K:=\Q(z)$ or in
its quadratic subfield $L$.  In particular, with $R$ [$R'$] as the set of
quadratic [non]residues mod $\ell$, if $z$ is a root of $f(x)$, then
\[ f(x) = \prod_{r \in R}(x-z^r);\ \ \ \ \ \ g(x)=\prod_{s \in R'}(x-z^s).  \]
\noindent 
Particular quantities are the {\bf trace-coefficients} of $f(x)$ and $g(x)$:
\be \label{eta}
   \eta := \sum_{r \in R}z^r;\ \ \ \ \ \ \eta' := \sum_{s \in R'}z^s.  
\eeq
That is, with $k := (\ell - 1)/2$ for ease of writing, the coefficients of
$x^{k-1}$ in $f(x)$ and $g(x)$, respectively, are $-\eta$ and $-\eta'$.

Also, $L = \Q(\eta)$, and $1 + \eta +\eta' = 0$.  In fact, 

\[  (x-\eta)(x - \eta') = x^2 +x + (1 - \leg{-1}\ell)/4.  \]
\noindent
Note also that
\[ f^*(x) = \left\{ \begin{array}{ll}
                      -g(x)  & \mbox{$\ell = 4N - 1$} \\  
                       f(x)  & \mbox{$\ell = 4N + 1$}
                    \end{array}.
            \right.  \]
Of course the same holds with $f$ and $g$ interchanged.

\smallskip

We focus now on the $[\ell, k]$ code $A$ defined as that generated by 
$(x-1)f(x)$.  $A^{+}$ is generated by $f(x)$.  (And $B^+$ is generated by
$g(x)$.  Orthogonality is laid out in {\sc nfd}.\footnote{The codes $A_{\infty}$ and $B_{\infty}$ over the quadratic number-field $L$ are conjugates of each other.})  $A^+$ is mapped to
$A_{\infty} \subset L^{\ell +1}$ by the rule $a \longrightarrow a;a_{\infty}$ 
in which
\be \label{gamma}
  a_{\infty} := \gamma\sum_i a_i,   
\eeq
and $\gamma$ satisfies 
\[ \ell \gamma^2 = \leg{-1};   \]
we use the Legendre symbol here.  The purpose of this clarification is to
determine the sign on $\gamma.$   We'll prove

\bigskip

{\bf Proposition}.  {\it Define $\eta$ and $\eta'$ as in (\ref{eta}).  The code $A^+$ generated by $f(x)$ extends to $A_{\infty}$ with the use of $\gamma \in L$ as specified in (\ref{gamma}).  Then}
\[  \ell \gamma =  -(\eta - \eta').   \]

We emphasize that $\eta$ is the trace-coefficient of $f(x)$.

This result nails down the too-vague discussion on page 132 of {\sc nfd}, 
which was based on hazy assumptions.  The proof below will be
followed by a clarification of the proof of the Gleason-Prange
theorem.

\smallskip

\pr  Since $x^{\ell} - 1 = (x-1)f(x)g(x)$, we have
\be \label{zroot}
A := ((x-1)f(x)) = \left( \frac{x^{\ell}-1}{g(x)} \right) = 
\left\{\begin{array}{ll}
       \left( \frac{x^{\ell}-1}{f^*(x)} \right)  & \mbox{$\ell = 4N - 1$} \\  
       \left( \frac{x^{\ell}-1}{g^*(x)} \right)  & \mbox{$\ell = 4N + 1$}
      \end{array}.
\right.
\eeq

We first recall that $\sigma$ (\cite{NFD}, p. 131)
is this monomial transformation of $L^{\ell + 1}$: with $\ep_i := \leg{i}$
for $0 < i < \ell$ and $\ep_0, \ep_{\infty} \in \{1, -1\}$ to be determined,  
\[ (a_0, \ldots , a_i, \ldots ;a_{\infty})\sigma := (\ep_0a_{\infty}, \ldots, \ep_ia_{-1/i}, \ldots ;\ep_{\infty}a_0).    \]

We begin  on page 131 of {\sc nfd} at 
``{\sc Case 1}: $\ell \equiv -1 \pmod{4}$.''

The proof of Case 1 has two parts.  Part I is the proof that 
$\langle1,0\rangle\sigma \in A_{\infty}$. With this Part we'll prove the Proposition.

Part II is the proof that $A\sigma \subset  A_{\infty}$.  We'll redo it
below to simplify the proof of the Gleason-Prange theorem.

\smallskip

Part I is OK until the top of page 132, where $z$ is not clearly specified.
For $\langle0,c\rangle = (T(cz^i))_{0 \leq i < \ell}$ to be in $A$, which we have 
defined to be $((x-1)f(x))$, 
(\ref{zroot}) tells us that $z$ must be a root of $f(x)$ (not of $g_1(x)$ as
stated in {\sc nfd}).  (If $f(z)=0$ is not clear, see Case 2 just below.)  
Still, the equation (3) of  {\sc nfd} is correct with our
present definitions of $\eta$ and $\eta'$.  Thus, when $\ell = 4N - 1$,
and if we take $\epsilon = 1$,
\[ \ell \gamma = -(\eta - \eta').   \]  

\smallskip

We now take up Part I of ``{\sc Case 2}: $\ell \equiv +1 \pmod{4}$.''
We take $\ep = 1.$

On page 133 of {\sc nfd}, it may help to note that $\langle1,0\rangle$ 
differs from $A_{\infty}$ to $B_{\infty}$, i.e.,
\[ \langle1,0\rangle_A = (1, 1, \ldots , 1; \ell \gamma)   \]
and
\[ \langle1,0\rangle_B = (1, 1, \ldots , 1; -\ell \gamma).   \]
  
As before, the two subparts of Part I are: Ia, prove that 
$\langle1,0\rangle_A\sigma$
is orthogonal to $\langle1,0\rangle_B$;  and  Ib, prove that 
$\langle1,0\rangle_A\sigma$ is orthogonal to $B$. 

For Ia:  The change of signs in the infinite coordinate makes this happen.

For Ib: Now $A^{\perp} = B^+ = (g(x)).$ 
We set up 
\[  F_z(c) := (T(c), T(cz), \ldots , T(cz^{\ell-1});0)  \]
as a general element of $B$. Since $B = ((x^{\ell} - 1)/f(x))$, and 
$f^*(x) = f(x)$, we see from (\ref{zroot}) that $z$ is a root of $f(x)$.

 Since 
\[ \langle1,0\rangle_A\sigma = (\ell \gamma, \ldots , \leg{i}, \ldots; 1),  \]
the dot product $F_z(c)\cdot \langle 1,0 \rangle \sigma$ is
\[  T(c)(\ell \gamma + \sum_{r \in R}z^r - \sum_{s \in R'}z^s).   \]
It is the same as it was in the prior case.  In other words, under the setup 
here, 
$\ell \gamma +\eta - \eta' = 0$ if and only if 
$\langle1,0\rangle_A\sigma$ is orthogonal 
to $B_{\infty}.$  This proves the Proposition.                \qed 

\bigskip

Now we go to Parts II of the two cases, to clarify the proof of the 
Gleason-Prange theorem, the burden of which is to prove that $\sigma$
is an invariance of the codes.

{\sc Case 1}: $\ell \equiv -1 \pmod{4}$.  It remains to prove that
$A\sigma \subset A_{\infty}$.  
 For a general element of $A$ we take
\be \label{-1}
 \langle0,c\rangle\  :=   (T(c), T(cz), \ldots , T(cz^{\ell-1});0),       
\eeq
where $z$ is a root of $f(x)$, as we saw earlier.  We apply $\sigma$, as in
{\sc nfd}, p. 133, to get
\be \label{-2}  
  \langle0,c\rangle\sigma = (a_0, \ldots ,a_{\ell - 1};a_{\infty}).
\eeq
In {\sc nfd} we first verified that $a_{\infty}$ is correct.  (See the equation
just above the line beginning ``from (2) and (1).'')   That calculation 
benefits from the Proposition just proved and states a correct result.

The only remaining hurdle is to show that as a polynomial this vector 
is a multiple of $f(x)$.  Since  
\be  \label{-3}
     a_i = \epsilon_i T(cz^{-1/i}),  
\eeq
we defined (noting $a_0 = 0$) for all $c \in K$, 
\be  \label{-4}
        D(c) := \sum_{1 \leq i < \ell}\epsilon_iT(cz^{-1/i})z^i.  
\eeq
Our object now is to prove $D(c)$ is always 0.  

We note that $D$ is linear from $K$ to itself.  Departing from {\sc nfd}'s
use of the ``quadratic-residue'' invariance $\tau$, we therefore prove it 0 
for $c = z, z^2, \ldots, z^{\ell-1},$ a spanning set for $K/L$.      Thus
\begin{eqnarray*}
  D(z^j) & = & \sum_{1 \leq i < \ell}\epsilon_iT(z^{j-1/i})z^i   \\
         & = & \sum_{i}\epsilon_i \sum_{r \in R}z^{r(j-1/i)}z^i   \\
         & = & \sum_{i,r} z^{rj -r/i + i}.
\end{eqnarray*}          

The rest of the proof is the same as in {\sc nfd}, except that the little 
polynomial is now  $x^2 + (rj-k)x -r$.  Also, for each $k,r$ this polynomial 
has two distinct, or no, roots in $GF(\ell).$  If two roots, one is in $R$ 
and the other in $R'$.  And the proof does not need
that $f(z) = 0$.  This settles Case 1.

\medskip

{\sc Case 2}: $\ell \equiv +1 \pmod{4}.$ To prove: $A\sigma \subset A_{\infty}.$

We imitate Case 1, making the necessary changes.

\smallskip

We define the general element of $A$ as before, except that now $z$ must
be a root of $g(x)$:

\[  \langle0,c\rangle\  :=   (T(c), T(cz), \ldots , T(cz^{\ell-1});0).   \]
 
And, as before,

\[  v:=\  \langle0,c\rangle\sigma = (0, \ldots, \epsilon_iT(cz^{-1/i}), \ldots; T(c)).   \]

Now $v$ is in $A_{\infty}$ if and only if its finite part, as 
a polynomial, is a multiple of $f(x)$.  The roots of $f(x)$ are $z^s$ for
$s \in R'$.  So: fix $s \in R$ and see that $v \in A_{\infty}$ iff 
\[ D'(c) := \sum_{1 \leq i < \ell}\epsilon_iT(cz^{-1/i})z^{si} = 0    \]
for all $c \in K$. 
 
Proceeding just as before, we write, for $0 < j < \ell$,
\[   D'(z^j) = \sum_i \epsilon_i\sum_{r \in R}z^{rj -r/i +si}.   \]
This expression is a polynomial in $z$ in which the coefficient of $z^k$
is the sum of the $\epsilon_i$ for which there are $r \in R$ such that
\[   rj -r/i +si = k.   \]
This comes down to $si^2 +(rj-k)i -r = 0$, i.e., 
\[   i^2 +s^{-1}(rj-k)i - s^{-1}r = 0.   \]
Again, the constant term is in $R'$, so not only are there no double roots
for $i$, but also if it has roots, one is in $R$ and one is in $R'$. 
Thus $D'(z^j) = 0.$                                       \qed

\medskip

{\sc A loose end.}  ``[W]e are free to choose $\epsilon = 1$ or
$\epsilon = -1$'' (page 132).  Recall that  $\epsilon := \epsilon_0$ and
$\epsilon_{\infty} := \leg{-1}\epsilon_0$.   We have always taken $\ep = 1$. 
But we are {\it not} free to choose it to be $-1$, because we know that 
\[ \langle1,0\rangle \sigma = (\epsilon \ell \gamma, \ldots \epsilon_i, \ldots;\leg{-1}\epsilon)  \] 
is in the code $A_{\infty}$ when $\epsilon = 1$, meaning that the polynomial
of its finite part is a multiple of $f(x)$.  If it were also in $A_{\infty}$
when we chose $\epsilon = -1$, the corresponding polynomial would be the 
same as the first one but with the constant term, not equal 0, 
of opposite sign.   It could not also be a multiple of $f(x)$.

This concludes my comments on the proof of the Gleason-Prange theorem.

\bigskip

{\sc A red herring.}  In a cyclic code of length $n$, for any polynomial
$f(x)$ dividing $x^n - 1$, we say that a polynomial $a(x)$ of degree less 
than $n$ is {\it recursive for $f(x)$} iff 
\be \label{recur}
            a(x)f(x) \equiv 0 \pmod{x^n - 1}.
\eeq  
Thus the set of all such $a(x)$ is the cyclic code with generator polynomial
$(x^n - 1)/f(x)$.  

The definition of ``recursive'' in {\sc nfd} (p. 125) used the reverse of $f(x)$ in 
(\ref{recur}).  Later (p. 129) we wrote of code ``$A$, recursive for $g_1(x)$,
generated as ideal by $g_2(x)$'', so we had slipped into the definition in
(\ref{recur}).  

Except for possibly confusing readers, no damage was done, because we never 
used recursion in working with the codes.  The best use of it, as on page 129,
is that saying a code is ``recursive for $f(x)$'' can be more economical,
even clearer, than saying it is generated by $(x^n - 1)/f(x)$.

For example, consider the [31,5] binary cyclic codes.  They are recursive for
the irreducible polynomials of degree 5 over $\F_2$.  One such is
$1+x^2 +x^5.$ It is easier to understand ``cyclic code recursive for
$1+x^2 +x^5$'' than ``cyclic code generated by $(x^{31} - 1)/(1+x^2 +x^5)$,''
however one might express the latter polynomial.

\bigskip

{\sc Conclusions.}  This note clarifies some murky points and corrects some
misstatements in, and simplifies, the proof of the Gleason-Prange theorem 
in {\sc nfd}.  One of these points is the sign of 
$\gamma$ used in the extension of the codes to length $\ell + 1$.

We emphasize that our basic datum is the factorization of $x^{\ell} - 1$ as
$(x-1)f(x)g(x)$ defined earlier.  This factorization determines the values
$\eta$ and $\eta'$, the trace-coefficients of the ``quadratic-residue''
polynomials $f(x)$ and $g(x)$, respectively.  These values in turn lead us
to the sign of $\gamma$.

\bigskip

{\sc Acknowledgment}.  The author thanks George Spelvin for calling his
attention to the warts in {\sc nfd} prettied up here.

\end{document}